\documentclass[11pt]{article}
\usepackage{epic,latexsym,amssymb}
\usepackage{amsfonts}
\usepackage{amscd}
\usepackage{amsmath}
\usepackage{graphicx}
\usepackage{color}
\usepackage{caption,subcaption}
\usepackage{tikz}

\usepackage{graphicx}

\usepackage{float}

\let\oldenumerate\enumerate
\renewcommand{\enumerate}{
  \oldenumerate
  \setlength{\itemsep}{0pt}
  \setlength{\parskip}{0.75pt}
  \setlength{\parsep}{0pt}
}

\newtheorem{thm}{Theorem}
\newtheorem{conj}{Conjecture}
\newtheorem{ob}[thm]{Observation}
\newtheorem{prop}[thm]{Proposition}

\newtheorem{claim}{Claim}
\newtheorem{subclaim}{Claim}[claim]

\newcommand{\od}{{\rm od}}

\newcommand{\barS}{\overline{S}}

\newcommand{\proof}{\noindent\textbf{Proof. }}

\newcommand{\smallqed}{{\tiny ($\Box$)}}
\newcommand{\qed}{$\Box$}

\newenvironment{unnumbered}[1]{\trivlist
\item [\hskip \labelsep {\bf #1}]\ignorespaces\it}{\endtrivlist}

\newcommand{\QEDmark}{\mbox{\textsc{qed}}}
\newcommand{\proofStarter}[1]{\textsc{#1} }

\def\vertex(#1){\put(#1){\circle*{2}}}
\def\vertexo(#1){\put(#1){\circle{2}}}
\def\vert(#1){\put(#1){\circle*{1.5}}}
\def\verto(#1){\put(#1){\circle{1.5}}}
\def\lab(#1)#2{\put(#1){\makebox(0,0)[c]{#2}}}
\definecolor{DarkGreen}{rgb}{0.2, 0.6, 0.3}

\newcommand{\Mike}[1]{{\color{blue}#1}}

\begin{document}

\title{The Forcing Number of Graphs with Given Girth}
\author{$^{1,2}$Randy Davila and $^1$Michael Henning\thanks{Research supported in part by the South African National Research Foundation and the University of Johannesburg} \\
\\
$^1$Department of Pure and Applied Mathematics\\
University of Johannesburg \\
Auckland Park 2006, South Africa \\
\small {\tt Email: mahenning@uj.ac.za} \\
\\
$^2$Department of Mathematics \\
Texas State University \\
San Marcos, TX 78666, USA \\
\small {\tt Email: rrd32@txstate.edu}
}

\date{}
\maketitle

\begin{abstract}
In this paper, we study a dynamic coloring of the vertices of a graph $G$ that starts with an initial subset $S$ of colored vertices, with all remaining vertices being non-colored. At each discrete time interval, a colored vertex with exactly one non-colored neighbor forces this non-colored neighbor to be colored. The initial set $S$ is called a forcing set of $G$ if, by iteratively applying the forcing process, every vertex in $G$ becomes colored. The forcing number, originally known as the \emph{zero forcing number}, and denoted $F(G)$, of $G$ is the cardinality of a smallest forcing set of $G$. We study lower bounds on the forcing number in terms of its minimum degree and girth, where the girth $g$ of a graph is the length of a shortest cycle in the graph. Let $G$ be a graph with minimum degree $\delta \ge 2$ and girth~$g \ge 3$. Davila and Kenter [Theory and Applications of Graphs, Volume 2, Issue 2, Article 1, 2015] conjecture that $F(G) \ge \delta + (\delta-2)(g-3)$. This conjecture has recently been proven for $g \le 6$. The conjecture is also proven when the girth $g \ge 7$ and the minimum degree is sufficiently large. In particular, it holds when $g = 7$ and $\delta \ge 481$, when $g = 8$ and $\delta \ge 649$, when $g = 9$ and $\delta \ge 30$, and when $g = 10$ and $\delta \ge 34$. In this paper, we prove the conjecture for $g \in \{7,8,9,10\}$ and for all values of $\delta \ge 2$.
\end{abstract}


\vskip 1 cm
\section{Introduction}

Graph \emph{dynamic colorings} are graph colorings that may \emph{change} with respect to discrete time intervals. One of the most prominent dynamic coloring is the result of the \emph{forcing process} (originally called the \emph{zero forcing process}), and its associated graph invariant, the \emph{forcing number} (originally called the \emph{zero forcing number}). These concepts first appeared during a workshop on linear algebra in relation to the \emph{minimum rank problem} \cite{AIM-Workshop}, and since then have been related to \emph{domination} and \emph{independence} \cite{k-Forcing}, \emph{network infection} \cite{quantum1}, and \emph{complexity} \cite{zf_np2}, to name a few. We highlight that computing the forcing number for a general graph is $NP$-hard \cite{zf_np}, and as such, finding computationally efficient bounds in terms of easily computable graph properties is of particular interest; see, for example, \cite{k-Forcing, Dynamic Forcing, Genter1, Genter2}.

Throughout this paper all graphs will be consider simple, undirected, and finite. Let $G=(V,E)$ be a graph with order $n=|V(G)|$, and size $m=|E(G)|$. Let $v$ be a vertex in $G$. A \emph{neighbor} of $v$ is a vertex adjacent to $v$. The \emph{open neighborhood} of $v$ in $G$, denoted $N_G(v)$, is the set of all neighbors of $v$ in $G$, whereas the \emph{closed neighborhood} of $v$ in $G$ is $N_G[v] = N_G(v)\cup \{v\}$. If the graph $G$ is clear from the context, we simply write $N(v)$ and $N[v]$ rather than $N_G(v)$ and $N_G[v]$, respectively. For a set $S \subseteq V$, its \emph{open neighborhood} is the set $N_G(S) = \bigcup_{v \in S} N(v)$, and its \emph{closed neighborhood} is the set $N_G[S] = N_G(S) \cup S$.

We denote the \emph{degree} of a vertex $v$ in $G$ by $d_G(v) = |N_G(v)|$. The minimum and maximum degrees among all vertices of $G$ are denoted by $\delta = \delta(G)$ and $\Delta = \Delta(G)$, respectively. For a subset $S$ of vertices in $G$, the subgraph induced by $S$ is denoted by $G[S]$. The degree of $v$ in $S$, denoted by $d_S(v)$, is the number of neighbors of $v$ in $G$ that belong to $S$. In particular, if $S = V(G)$, then $d_S(v) = d_G(v)$.
The \emph{distance} between two vertices $u$ and $v$ in $G$, denoted $d_G(u,v)$ or simply $d(u,v)$ if the graph $G$ is clear from context, is the minimum length of a $(u,v)$-path in $G$. The length of a shortest cycle in $G$ is the \emph{girth} of $G$, denoted by $g = g(G)$. We will use the notation $P_n$, $C_n$, $K_n$, and $K_{n,m}$, to denote the \emph{path} on $n$ vertices, the \emph{cycle} on $n$ vertices, the \emph{complete graph} on $n$ vertices, and the complete bipartite graph with partite sets of sizes $n$ and $m$, respectively. If $G$ does not contain a graph $F$ as an induced subgraph, we say that $G$ is $F$-\emph{free}. A graph is \emph{triangle}-\emph{free} if it is $K_3$-free.

Let $D$ be a digraph with vertex set $V(D)$ and arc set $A(D)$. The \emph{out-degree} of a vertex $v$ in $D$, denoted by $\od_D(v)$, is the number of vertices $w$ such that $(v,w) \in A(D)$, where $(v,w)$ denotes an arc from $v$ to $w$; that is, $\od_D(v) = |\{ w \in V(D) \colon (v,w) \in A(D)\}|$. We denote the number of arc in a digraph $D$ by $m(D)$, and note that $m(D) = \sum_{v \in V(D)} \od_D(v)$.

Adopting the notation of~\cite{DaHeMaPe16}, given a graph $G$, the \emph{forcing process} is defined as follows: Let $S \subseteq V(G)$ be an initial set of ``colored'' vertices; all remaining vertices being ``non-colored''. A vertex in a set $S$, we call $S$-\emph{colored}, while a vertex not in $S$ we call $S$-\emph{uncolored}. At each time step, a colored vertex $v$ with exactly one non-colored neighbor will change, or \emph{force}, the non-colored neighbor to be colored. We call such a vertex $v$ a \emph{forcing colored vertex}, or simply a \emph{forcing vertex}. Further, at the time when the vertex $v$ forces its non-colored neighbor to be colored, we say that the vertex $v$ is \emph{played}. A set $S \subseteq V(G)$ of initially colored vertices is called a \emph{forcing set} if, by iteratively applying the forcing process, all of $V(G)$ becomes colored. We call such a set $S$ an  $S$-\emph{forcing set}. The \emph{forcing number} of a graph $G$, denoted by $F(G)$, is the cardinality of a smallest forcing set. If $S$ is a forcing set in $G$ and $v$ is an $S$-colored vertex that forces a new vertex to be colored, then we call $v$ an $S$-\emph{forcing vertex}.

\section{Main Result}

In this paper, we study the following intriguing conjecture posed by Davila and Kenter~\cite{Davila Kenter}.

\begin{conj}{\rm (\cite{Davila Kenter})}
\label{conj2}
If $G$ is a graph with girth $g \ge 3$ and minimum degree $\delta\ge 2$, then
\begin{equation*}
F(G) \ge \delta + (\delta-2)(g-3).
\end{equation*}
\end{conj}

Gentner, Penso, Rautenbach, and Souzab~\cite{Genter1} and Gentner and Rautenbach~\cite{Genter2} have shown that Conjecture~\ref{conj2} is true for small girth $g \le 6$, while Davila and Kenter~\cite{Davila Kenter} have proven that Conjecture~\ref{conj2} is true for girth $g \ge 7$ and sufficiently large minimum degree. We state these results formally as follows.

\begin{thm}\label{c:Girth7+}
If $G$ is a graph with girth $g \ge 3$ and minimum degree~$\delta \ge 2$, then the following holds.
\\
\indent {\rm (a)} {\rm (\cite{Genter1,Genter2})} If $g \le 6$, then Conjecture~\ref{conj2} is true.
\\
\indent {\rm (b)} {\rm (\cite{Davila Kenter})} If $g = 7$ and $\delta \ge 481$, then Conjecture~\ref{conj2} is true.
\\
\indent {\rm (c)} {\rm (\cite{Davila Kenter})} If $g = 8$ and $\delta \ge 649$, then Conjecture~\ref{conj2} is true. \\
\indent {\rm (d)} {\rm (\cite{Davila Kenter})} If $g = 9$ and $\delta \ge 30$, then Conjecture~\ref{conj2} is true.
\\
\indent {\rm (e)} {\rm (\cite{Davila Kenter})} If $g = 10$ and $\delta \ge 34$, then Conjecture~\ref{conj2} is true.
\end{thm}

Our aim in this paper is to prove that Conjecture~\ref{conj2} is true when the girth $g \in \{7,8,9,10\}$, for all minimum degree~$\delta \ge 2$. This improves the result of Theorem~\ref{c:Girth7+}(b),~\ref{c:Girth7+}(c) and~\ref{c:Girth7+}(d) which imposes a restriction on the minimum degree~$\delta$. We state our result formally as follows.

\begin{thm}
\label{thm:main1a}
If $G$ is a graph with girth $g \in \{7,8,9,10\}$ and minimum degree~$\delta \ge 2$, then Conjecture~\ref{conj2} is true.
\end{thm}

\section{Known Results and Motivation}

As remarked earlier, finding bounds on $F(G)$ in terms of easily computable graph properties is of interest. The earliest such bound is given in the original paper~\cite{AIM-Workshop} which showed that $F(G)$ is at least as small as the minimum degree. In particular, we note that any initially forcing vertex must be colored along with all but one of its neighbors. We state this result formally with the following proposition.

\begin{prop}{\rm (\cite{AIM-Workshop})}
\label{delta}
If $G$ is a graph with minimum degree $\delta$, then $F(G)\ge \delta(G)$, and this bound is sharp.
\end{prop}

Proposition~\ref{delta} is sharp, as can be seen by considering the path $P_n$, the cycle $C_n$, and the complete graph $K_n$. Since the leaf of every non-trivial path is a forcing set, we note that $F(P_n) = 1$. As observed in~\cite{DaHeMaPe16}, paths are the only graphs $G$ satisfying $F(G) = 1$, and are therefore a special class of graphs when considering $F(G)$. Complete graphs have the largest possible forcing number. As observed in~\cite{DaHeMaPe16}, if $G$ is a connected graph of order~$n \ge 2$, then $F(G) \le n-1$, with equality if and only if $G = K_n$. We remark that this upper bound of~$n-1$ implies that there exists at least one played vertex in any minimum forcing set for non-empty graphs.

Barioli et al.~\cite{Barioli13} prove that the forcing number of a graph is at least its tree-width. Using this result, and a result of Chandrana and Subramanian~\cite{ChSu05} that establishes a lower bound on the tree-width of a graph in terms of its average degree and girth, Davila and Kenter~\cite{Davila Kenter} prove that Conjecture~\ref{conj2} is true for graphs with girth at least~$7$ and sufficiently large minimum degree. More precisely, they prove the following result.

\begin{thm}{\rm (\cite{Davila Kenter})} \label{t:Girth7+}
Conjecture~\ref{conj2} is true for all graphs with minimum degree~$\delta$ and given girth~$g$, where $g \ge 7$, that satisfy
\[
\frac{(\delta - 1)^{\lfloor (g-1)/2 \rfloor - 1}}{12(g+1)} \ge \delta + (\delta - 2)(g-3).
\]
\end{thm}

We remark that Theorem~\ref{c:Girth7+}(b),~\ref{c:Girth7+}(c), ~\ref{c:Girth7+}(d) and~\ref{c:Girth7+}(e) follows from Theorem~\ref{t:Girth7+} in the special case when $g \in \{7,8,9,10\}$.

\section{Proof of Main Result}

Let $G$ be a graph with minimum degree~$\delta \ge 2$ and girth $g$, where $g \ge 5$, and consider a shortest cycle $C$ in $G$ of length~$g$. The girth at least~$5$ constraint implies that no two vertices on $C$ have a common neighbor outside $C$, implying that $G$ has at least $g(\delta - 1)$ vertices. We state this observation formally as follows.

\begin{ob}
\label{o:girth}
If $G$ is a graph of order~$n$ with minimum degree~$\delta \ge 2$ and girth $g$, where $g \ge 5$, then $n \ge g (\delta - 1)$.
\end{ob}

In order to proceed with the proof of Theorem~\ref{thm:main1a}, we recall a classical result known as \emph{Mantel's Theorem}.

\begin{thm}[Mantel's Theorem]
If a graph $G$ on $n$ vertices is triangle-free, then it contains at most $\frac{n^2}{4}$ edges.
\end{thm}

We are now ready to give the proof of Theorem~\ref{thm:main1a}.

\noindent \textbf{Theorem~\ref{thm:main1a}} \emph{If $G$ is a graph with girth $g \in \{7,8,9,10\}$ and minimum degree~$\delta \ge 2$, then Conjecture~\ref{conj2} is true.}

\noindent
\textbf{Proof of Theorem~\ref{thm:main1a}.}
Let $G$ be a graph of order~$n$ with girth $g$, where $g \in \{7,8,9,10\}$, and minimum degree $\delta \ge 2$. We wish to show that $F(G) \ge \delta + (\delta-2)(g-3)$. Suppose, to the contrary, that $F(G) \le \delta + (\delta-2)(g-3) - 1$. If $\delta = 2$, then $\delta + (\delta-2)(g-3) = \delta$, and so, by Proposition~\ref{delta}, $F(G) \ge \delta + (\delta-2)(g-3)$, a contradiction. Hence, $\delta \ge 3$. By Observation~\ref{o:girth}, $n \ge g(\delta - 1)$.
Let $S \subseteq V(G)$ be a minimum forcing set of $G$, and so $|S| = F(G) \le \delta + (\delta-2)(g-3) - 1$. We state our supposition formally as follows.

\begin{claim}\label{c:1}
The following holds. \\[-28pt]
\begin{enumerate}
\item If $g = 7$, then $|S| \le 5\delta - 9$.
\item If $g = 8$, then $|S| \le 6\delta - 11$.
\item If $g = 9$, then $|S| \le 7\delta - 13$.
\item If $g = 10$, then $|S| \le 8\delta - 15$.
\end{enumerate}
\end{claim}

Let $\barS = V(G) \setminus S$, and so $\barS$ is the set of all $S$-uncolored vertices. Thus, $|\barS| =  n - |S| \ge g(\delta - 1) - \delta - (\delta-2)(g-3) + 1  = g + (2\delta - 5) \ge g + 1$. Since $S$ is a forcing set of $G$, there is a sequence $x_1,x_2,\ldots,x_t$ of played vertices in the forcing process that results in all $V(G)$ colored, where $x_i$ denotes the forcing colored vertex played in the $i$th step of the process. We note that $t = |\barS| > g - 2$. Let $X = \{x_1,x_2,\ldots,x_{g-2}\}$. For $i \in [g-2]$, let $X_i = N[x_i]$ and let
\[
X_{\le i} = \bigcup_{j=1}^{i} X_j.
\]

Let $S_1$ be the neighbors of $x_1$ in $S$, let $S_2$ be the neighbors of $x_2$ in $S$ that do not belong to $X_1$, let $S_3$ be the neighbors of $x_3$ in $S$ that do not belong to $X_1 \cup X_2$, and so on. Thus, $S_1 = S \cap N(x_1)$, and for $i \in [g-2] \setminus \{1\}$, the set $S_i$ is the set of neighbors of $x_i$ in $S$ that do not belong to $X_{\le i-1}$; that is,
\[
S_i = S \cap (N(x_i) \setminus X_{\le i-1} ).
\]

We note that $S_i \subset X_i \setminus \{x_i\}$ for each $i \in [g-2]$. Let
\[
S_X^* = \bigcup_{i=1}^{g-2} S_i \hspace*{0.75cm} \mbox{and} \hspace*{0.75cm} S_X = X \cap (S \setminus S_X^*).
\]

By definition, $S_i \cap S_j = \emptyset$ for $1 \le i,j \le [g-2]$ and $i \ne j$, and so
\[
|S_X^*| = \sum_{i=1}^{g-2} |S_i|.
\]
\indent
Since the girth of $G$ is greater than~$4$, the vertex $x_i$ has at most one neighbor in $X_j$ for each $j < i$ and $i \in [g-2]$. Let $D$ be the digraph with vertex set $V(D) = X$ and with arc set $A(D)$ defined as follows. For each $i \in [g-2]$, we add an arc from $x_i$ to $x_j$ if $i > j$ and $x_i$ has a neighbor in $X_j$ in the graph $G$.
We will now show a number of claims which culminate with a contradiction to $G$ being a counterexample.
\begin{claim}\label{c:2}
The following holds. \\[-28pt]
\begin{enumerate}
\item $|S| \ge |S_X| + |S_X^*|$.
\item $|S_X| \ge 1$ and $x_1 \in S_X$.
\item $|S_X^*| \ge (g-2)(\delta - 1) - m(D)$.
\end{enumerate}
\end{claim}
\proof Part~(a) follows immediately from the fact that the sets $S_X^*$ and $S_X$ are vertex disjoint and $S \cap N_G[X] = S_X \cup S_X^*$. By definition of the sets $S_i$, $i \in [g-2]$, the first vertex played, namely $x_1$, does not belong to $S_X^*$. However, the vertex $x_1$ is an $S$-forcing vertex and therefore belongs to $S$, and so $x_1 \in S_X$. This establishes Part~(b). To prove part~(c), we show that $|S_i| = d_G(x_i) - \od_{D}(x_i) - 1$ for all $i \in [g-2]$. Since the vertex $x_1$ has exactly one $S$-uncolored neighbor, this implies that $|S_1| = d_G(x_1) - \od_{D}(x_1) - 1$ noting that $\od_{D}(x_1) = 0$. Hence we may assume that $i \ge 2$. The vertex $x_i$ has exactly $\od_{D}(x_i)$ neighbors in $X_{\le i-1}$ (with at most one neighbor in $X_j$ for each $j \in [i-1]$).  Further, since $x_i$ is a forcing vertex in the $i$th step of the forcing process, it has exactly one neighbor in $\barS \setminus X_{\le i-1}$. Thus, $|S_i| = d_G(x_i) - \od_{D}(x_i) - 1$. Therefore,
\[
|S_X^*| = \sum_{i=1}^{g-2} |S_i| = \sum_{i=1}^{g-2} (d_G(x_i) - \od_{D}(x_i) - 1) \ge (g-2)(\delta - 1) - m(D),
\]
noting that $\displaystyle{\sum_{i=1}^{g-2} \od_{D}(x_i) = m(D)}$. This proves Part~(c).~\smallqed

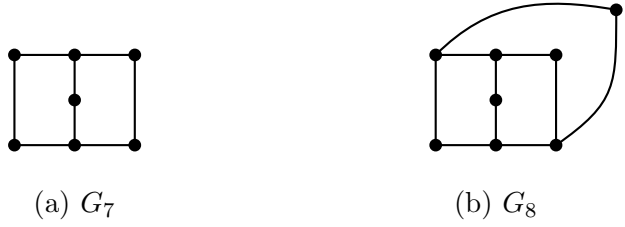
\begin{figure}[htb]
\begin{center}
\begin{tikzpicture}[scale=.8,style=thick,x=1cm,y=1cm]
\def\vr{2.5pt} 
\path (0,0) coordinate (a);
\path (1,0) coordinate (b);
\path (2,0) coordinate (c);
\path (2,1.5) coordinate (d);
\path (1,1.5) coordinate (e);
\path (0,1.5) coordinate (f);
\path (1,0.75) coordinate (g);
%
\draw (a) -- (b);
\draw (b) -- (c);
\draw (c) -- (d);
\draw (d) -- (e);
\draw (e) -- (f);
\draw (f) -- (a);
\draw (g) -- (b);
\draw (g) -- (e);
\draw (a) [fill=black] circle (\vr);
\draw (b) [fill=black] circle (\vr);
\draw (c) [fill=black] circle (\vr);
\draw (d) [fill=black] circle (\vr);
\draw (e) [fill=black] circle (\vr);
\draw (f) [fill=black] circle (\vr);
\draw (g) [fill=black] circle (\vr);
\draw (1,-1) node {(a) $G_7$};
\path (7,0) coordinate (a);
\path (8,0) coordinate (b);
\path (9,0) coordinate (c);
\path (9,1.5) coordinate (d);
\path (8,1.5) coordinate (e);
\path (7,1.5) coordinate (f);
\path (8,0.75) coordinate (g);
\path (10,2.25) coordinate (h);
%
\draw (a) -- (b);
\draw (b) -- (c);
\draw (c) -- (d);
\draw (d) -- (e);
\draw (e) -- (f);
\draw (f) -- (a);
\draw (g) -- (b);
\draw (g) -- (e);
\draw (a) [fill=black] circle (\vr);
\draw (b) [fill=black] circle (\vr);
\draw (c) [fill=black] circle (\vr);
\draw (d) [fill=black] circle (\vr);
\draw (e) [fill=black] circle (\vr);
\draw (f) [fill=black] circle (\vr);
\draw (g) [fill=black] circle (\vr);
\draw (h) [fill=black] circle (\vr);
\draw (h) to[out=-190,in=45, distance=1.25cm ] (f); 
\draw (h) to[out=270,in=35, distance=1.25cm ] (c); 
\draw (8,-1) node {(b) $G_8$};
\end{tikzpicture}
\end{center}
\vskip -0.4cm
\caption{The graphs $G_7$ and $G_8$.} \label{f:G7}
\end{figure}

\medskip
Let $G_D$ denote the underlying graph of the digraph $D$ obtained from $D$ by removing the directions of the arcs in $A(D)$. We note that $V(G_D) = X$ and $m(G_D) = m(D)$.
Further, we note that if $x_{i_1}x_{i_2}$ is an edge in $G_D$, then either $x_{i_1}$ and $x_{i_2}$ are adjacent in $G$ or $x_{i_1}$ and $x_{i_2}$ have a common neighbor in $G$. Let $G_7$ and $G_8$ be the graphs illustrated in Figure~\ref{f:G7}(a) and~\ref{f:G7}(b), respectively.

\begin{claim}\label{c:3}
The following holds. \\[-28pt]
\begin{enumerate}
\item If $g = 7$, then $|S_X^*| \ge 5\delta - 11$. Further, if  $|S_X^*| = 5\delta - 11$, then $G_D \cong K_{2,3}$.
\item If $g = 8$, then $|S_X^*| \ge 6\delta - 15$. Further, if $|S_X^*| = 6\delta - 15$, then $G_D \cong K_{3,3}$.
\item If $g = 9$, then $|S_X^*| \ge 7\delta - 15$. Further, if $|S_X^*| = 7\delta - 15$, then  $G_D \cong G_7$.
\item If $g = 10$, then $|S_X^*| \ge 8\delta - 18$. Further, if $|S_X^*| = 8\delta - 18$, then  $G_D \cong G_8$.
\end{enumerate}
\end{claim}
\proof Suppose first that $g\in \{7,8\}$. In this case, we note that the graph $G_D$ is triangle-free. Thus, by Mantel's Theorem, the maximum number of edges in the graph $G_D$, which has order~$g-2$, is~$\lfloor (g-2)^2/4 \rfloor$. Further, the graph $K_{\lfloor (g-2)/2 \rfloor , \lceil (g-2)/2 \rceil}$ is the unique extremal graph. This observation, together with Claim~\ref{c:2}(c) and the fact that $m(D) = m(G_D) \le \lfloor (g-2)^2/4 \rfloor$, yields Parts~(a) and~(b). In particular, we note that if $g = 7$, then $m(D) \le 6$, while if $g = 8$, then $m(D) \le 9$.

Suppose next that $g \in \{9,10\}$. In this case, we note that the graph $G_D$, which has order~$g-2$, contains neither three-cycles nor four-cycles. It is well-known (see, for example,~\cite{GaKwLa93}, or simply use a computer) that if $g = 9$, then such a graph $G_D$ has at most eight edges and there is a unique extremal graph of size~$8$, namely~$G_7$. Further, if $g = 10$, then such a graph $G_D$ has at most ten edges and there is a unique extremal graph of size~$10$, namely~$G_8$.~\smallqed

\begin{claim}\label{c:4}
For $i \in [g-2]$, the following holds. \\[-28pt]
\begin{enumerate}
\item If $x_i \in \barS$, then $x_i$ is adjacent in $G$ to $x_j$ for some $j \in [i-1]$.
\item If $x_i$ is adjacent in $G$ to no vertex $x_j$ where $j \in [i-1]$, then $x_i \in S_X$.
\end{enumerate}
\end{claim}
\proof
By definition of the forcing process, if $x_i \in \barS$ for some $i \in [g-2]$, then $i \ge 2$ and there exists some $j \in [i-1]$ such that the vertex $x_i$ is the vertex in $\barS$ that becomes colored when the vertex $x_j$ is played, implying that $x_i$ is adjacent to $x_j$ in $G$. This proves Part~(a). To prove Part~(b), suppose that $x_i$ is adjacent to no vertex $x_j$ where $j \in [i-1]$. By Part~(a), $x_i \in S$. Since $x_i$ is not adjacent to $x_j$ in $G$ for any $j \in [i-1]$, we note that $x_i \notin S_j$. This implies, by definition of the set $S_X^*$, that the vertex $x_i$ does not belong to $S_X^*$, implying that $x_i \in S_X$.~\smallqed

\newpage
\begin{claim}\label{c:5a}
$g \ne 9$.
\end{claim}
\proof Suppose, to the contrary, that $g = 9$. As shown in the proof of Claim~\ref{c:3}(c), the graph $G_D$, of order~$7$, contains neither three-cycles nor four-cycles and satisfies $m(G_D) \le 8$.

Suppose that $m(G_D) = 8$, implying that $G_D \cong G_7$. By Claim~\ref{c:3}(c), $|S_X^*| \ge 7\delta - 15$. Thus, by Claim~\ref{c:1} and Claim~\ref{c:2}(a), $7\delta - 13 \ge |S| \ge |S_X| + |S_X^*| \ge |S_X| + 7\delta - 15$, implying that $|S_X| \le 2$. By Claim~\ref{c:2}(b), $x_1 \in S_X$. Thus for at least five values of $i \in [7] \setminus \{1\}$, the vertex $x_i$ does not belong to $S_X$, and is therefore by Claim~\ref{c:4} adjacent to some vertex $x_j$ where $j \in [i-1]$. This implies that at least five edges in $G_D$ correspond to edges in the graph $G$. Recall that every edge $x_{i_1}x_{i_2}$ in $G_D$ is an edge in $G$ or corresponds to a path $x_{i_1}vx_{i_2}$ where $v$ is a common neighbor of $x_{i_1}$ and $x_{i_2}$ in $G$. Therefore, as least one of the two $5$-cycles in $G_D$ corresponds to a cycle in $G$ of length less than~$9$, a contradiction to our supposition that $g = 9$. Hence, $m(G_D) \le 7$.

Since $m(D) = m(G_D) \le 7$, Claim~\ref{c:2}(c) implies that $|S_X^*| \ge 7\delta - 14$. Thus, by Claim~\ref{c:1} and Claim~\ref{c:2}(c), $7\delta - 13 \ge |S| \ge |S_X| + |S_X^*| \ge 7\delta - 13$. Hence, we must have equality throughout this inequality chain, implying that $|S| = 7\delta - 13$, and so $|S_X| = 1$ and $|S_X^*| = 7\delta - 14$. Further, $m(D) = m(G_D) = 7$ and $S_X = \{x_1\}$. Thus for all $i \in [7] \setminus \{1\}$, the vertex $x_i \notin S_X$, implying that $x_i$ is adjacent to some vertex $x_j$ where $j \in [i-1]$. By the girth $g = 9$ supposition, this implies that each such vertex $x_i$ is adjacent to exactly one vertex $x_j$ where $j \in [i-1]$. Therefore, $G_D$ is a tree, and so $m(G_D) = 6$, a contradiction.~\smallqed

\begin{claim}\label{c:5b}
$g \ne 10$.
\end{claim}
\proof Suppose, to the contrary, that $g = 10$. As shown in the proof of Claim~\ref{c:3}(d), the graph $G_D$, of order~$8$, contains neither three-cycles nor four-cycles and satisfies $m(G_D) \le 10$.

Suppose that $m(G_D) = 10$, implying that $G_D \cong G_8$. By Claim~\ref{c:3}(d), $|S_X^*| \ge 8\delta - 18$. Thus, by Claim~\ref{c:1} and Claim~\ref{c:2}(a), $8\delta - 15 \ge |S| \ge |S_X| + |S_X^*| \ge |S_X| + 7\delta - 18$, implying that $|S_X| \le 3$. By Claim~\ref{c:2}(b), $x_1 \in S_X$. Thus for at least five values of $i \in [8] \setminus \{1\}$, the vertex $x_i$ does not belong to $S_X$, and is therefore by Claim~\ref{c:4} adjacent to some vertex $x_j$ where $j \in [i-1]$. This implies that at least five edges in $G_D$ correspond to edges in the graph $G$. Every $5$-cycle in $G_D$ corresponds to a cycle in $G$ of length at most~$10$. Thus, since the girth of $G$ is~$10$, every edge in a $5$-cycle in $G_D$ corresponds to a path of length~$2$ in $G$. However, every edge of $G_D$ belongs to some $5$-cycle, implying that no edge of $G_D$ corresponds to an edge of $G$ and therefore that $S_X = X$, and so $|S_X| = 8$, a contradiction.  Hence, $m(G_D) \le 9$. 

If $m(G_D) \le 7$, then Claim~\ref{c:2}(c) implies that $|S_X^*| \ge 8\delta - 15$. Thus, by Claim~\ref{c:2}, $|S| \ge |S_X| + |S_X^*| \ge 8\delta - 14$, contradicting Claim~\ref{c:1}(d). Hence, $m(G_D) \in \{8,9\}$, implying that $G_D$ contains a cycle. 

Since $m(D) = m(G_D) \le 9$, Claim~\ref{c:2}(c) implies that $|S_X^*| \ge 8\delta - 17$. Thus, by Claim~\ref{c:1} and Claim~\ref{c:2}(a), $8\delta - 15 \ge |S| \ge |S_X| + |S_X^*| \ge |S_X| + 8\delta - 17$, implying that $|S_X| \le 2$. Thus for at least six values of $i \in [8] \setminus \{1\}$, the vertex $x_i$ does not belong to $S_X$, and is therefore by Claim~\ref{c:4} adjacent to some vertex $x_j$ where $j \in [i-1]$. This implies that at least six edges in $G_D$ correspond to edges in the graph $G$. As observed earlier, $G_D$ contains a cycle, but contains neither three-cycles nor four-cycles. Further, $m(G_D) \in \{8,9\}$.

Suppose that $G_D$ contains a $5$-cycle. Since the girth of $G$ is~$10$, this implies that no edge of the $5$-cycle corresponds to an edge of $G$. Thus, there are at most~$m(D) - 5 \le 9 - 5 = 4$ edges of $G_D$ that correspond to edges in the graph $G$, a contradiction. Hence, $G_D$ has girth at least~$6$. Suppose that $G_D$ contains a $6$-cycle. In this case, since the girth of $G$ is~$10$ at least four edges of the $6$-cycle do not corresponds to an edge of $G$. Thus, there are at most~$m(D) - 4 \le 9 - 4 = 5$ edges of $G_D$ that correspond to edges in the graph $G$, a contradiction. Hence, $G_D$ has girth at least~$7$, implying that $m(G_D) = 8$. 

Claim~\ref{c:2}(c) implies that $|S_X^*| \ge 8\delta - 16$. Thus, by Claim~\ref{c:1} and Claim~\ref{c:2}(d), $8\delta - 15 \ge |S| \ge |S_X| + |S_X^*| \ge 8\delta - 15$. Hence, we must have equality throughout this inequality chain, implying that $|S| = 8\delta - 15$, and so $|S_X| = 1$ and $S_X = \{x_1\}$. Thus for all $i \in [8] \setminus \{1\}$, the vertex $x_i \notin S_X$, implying that $x_i$ is adjacent to some vertex $x_j$ where $j \in [i-1]$. By the girth $g = 10$ supposition, this implies that each such vertex $x_i$ is adjacent to exactly one vertex $x_j$ where $j \in [i-1]$. Therefore, $G_D$ is a tree, and so $m(G_D) = 7$, a contradiction.~\smallqed

\medskip
By Claims~\ref{c:5a} and~\ref{c:5b}, $g \in \{7,8\}$.

\begin{claim}\label{c:6}
The vertices $x_1$ and $x_2$ are not adjacent.
\end{claim}
\proof Suppose, to the contrary, that the vertices $x_1$ and $x_2$ are adjacent.

\begin{subclaim}\label{c:6.1}
The vertex $x_3$ is adjacent to neither $x_1$ nor $x_2$.
\end{subclaim}
\proof  Suppose, to the contrary, that $x_3$ is adjacent to $x_1$ or $x_2$. The girth at least seven requirement implies that the vertex $x_i$ has at most one neighbor in $X_{\le 3}$ for every $i \in \{4,\ldots,g-2\}$. Thus, if $g = 7$, then $m(D) \le 5$, while if $g = 8$, then, recalling that $G_D$ is triangle-free and therefore there are at most two edges in $G_D[\{x_4,x_5,x_6\}]$, this implies that  $m(D) \le 7$. Hence, by Claim~\ref{c:2}(c), if $g = 7$, then $|S_X^*| \ge 5\delta - 10$, while if $g = 8$, then $|S_X^*| \ge 6\delta - 13$.

\begin{unnumbered}{Claim~\ref{c:6.1}.1}
$g = 8$.
\end{unnumbered}
\proof
Suppose that $g = 7$. In this case, $|S_X^*| \ge 5\delta - 10$. By Claim~\ref{c:2}(b), $|S_X| \ge 1$. Thus, by Claim~\ref{c:1} and Claim~\ref{c:2}(a), $5\delta - 9 \ge |S| \ge |S_X| + |S_X^*| \ge 5\delta - 9$. Hence, we must have equality throughout this inequality chain, implying that $|S| = 5\delta - 9$, and so $|S_X| = 1$ and $|S_X^*| = 5\delta - 10$. Further, $S_X = \{x_1\}$. Since $x_4 \notin S_X$, the vertex $x_4$ is adjacent to some vertex $x_j$ where $j \in [3]$. By the girth condition, such a vertex $x_j$ is the only neighbor of $x_4$ in $X_{\le 3}$. Since $x_5 \notin S_X$, the vertex $x_5$ is adjacent to some vertex $x_r$ where $r \in [4]$. This implies by the girth seven requirement that $x_5$ has exactly one neighbor in $X_{\le 4}$, and therefore that $m(D) \le 4$. Hence, by Claim~\ref{c:2}(c), $|S_X^*| \ge 5\delta - 9$, a contradiction.~\smallqed

\medskip
By Claim~\ref{c:6.1}.1, the girth $g = 8$. Thus, $|S_X^*| \ge 6\delta - 13$.

\begin{unnumbered}{Claim~\ref{c:6.1}.2}
The vertex $x_4$ is adjacent to no vertex $x_j$ where $j \in [3]$.
\end{unnumbered}
\proof
Suppose, to the contrary, that $x_4$ is adjacent to some vertex $x_j$, where $j \in [3]$. Such a vertex $x_j$ is the only neighbor of $x_4$ in $X_{\le 3}$. The girth eight requirement implies that each of $x_5$ and $x_6$ has at most one neighbor in $X_{\le 4}$, and therefore that $m(D) \le 6$. Hence, by Claim~\ref{c:2}(c), $|S_X^*| \ge 6\delta - 12$. By Claim~\ref{c:2}(b), $|S_X| \ge 1$.  Thus, by Claim~\ref{c:1} and Claim~\ref{c:2}(a), $5\delta - 11 \ge |S| \ge |S_X| + |S_X^*| \ge 6\delta - 11$. Hence, we must have equality throughout this inequality chain, implying that $|S| = 6\delta - 11$, and so $|S_X| = 1$ and $|S_X^*| = 6\delta - 12$. Further, $S_X = \{x_1\}$. Since neither $x_5$ nor $x_6$ belongs to $S_X$, Claim~\ref{c:4}(b) implies that $x_i$ is adjacent to some vertex $x_j$ where $j \in [i-1]$ for $i \in \{5,6\}$. This, together with the girth eight requirement, implies that the vertex $x_6$ has exactly one neighbor in $X_{\le 5}$. This in turn implies that $m(D) \le 5$. Hence, by Claim~\ref{c:2}(c), $|S_X^*| \ge 6\delta - 11$, a contradiction.~\smallqed

\medskip
By Claim~\ref{c:6.1}.2, the vertex $x_4$ is adjacent to no vertex $x_j$ where $j \in [3]$. By Claim~\ref{c:4}(b), $x_4 \in S_X$, and so $|S_X| \ge 2$. As observed earlier, $|S_X^*| \ge 6\delta - 13$. Thus, by Claim~\ref{c:1} and Claim~\ref{c:2}(a), $5\delta - 11 \ge |S| \ge |S_X| + |S_X^*| \ge 6\delta - 11$. Hence, we must have equality throughout this inequality chain, implying that $|S| = 6\delta - 11$, and so $|S_X| = 2$ and $|S_X^*| = 6\delta - 13$. Further, $S_X = \{x_1,x_4\}$.

\begin{unnumbered}{Claim~\ref{c:6.1}.3}
The vertex $x_4$ has no neighbor in $X_{\le 3}$.
\end{unnumbered}
\proof
Suppose, to the contrary, that $x_4$ has a neighbor in $X_{\le 3}$. By the girth condition, this is the only neighbor of $x_4$ in $X_{\le 3}$. Since neither $x_5$ nor $x_6$ belongs to $S_X$, Claim~\ref{c:4}(b) implies that $x_i$ is adjacent to some vertex $x_j$ where $j \in [i-1]$ for $i \in \{5,6\}$. This, together with the girth eight requirement, implies that the vertex $x_5$ has exactly one neighbor in $X_{\le 4}$, and the vertex $x_6$ has exactly one neighbor in $X_{\le 5}$. Therefore, $m(D) \le 6$. Hence, by Claim~\ref{c:2}(c), $|S_X^*| \ge 6\delta - 12$, a contradiction.~\smallqed

\medskip
By Claim~\ref{c:6.1}.3, the vertex $x_4$ has no neighbor in $X_{\le 3}$, implying by our earlier observations that $m(D) \le 6$ and therefore that $|S_X^*| \ge 6\delta - 12$, a contradiction. This completes the proof of Claim~\ref{c:6.1}.~\smallqed

\medskip
By Claim~\ref{c:6.1}, the vertex $x_3$ is adjacent to neither $x_1$ nor $x_2$. Hence by Claim~\ref{c:4}(b), $x_3 \in S_X$, and so $|S_X| \ge 2$.

\begin{subclaim}\label{c:6.2}
The vertex $x_3$ has no neighbor in $X_{\le 2}$.
\end{subclaim}
\proof
Suppose, to the contrary, that $x_3$ has a neighbor in $X_{\le 2}$. Suppose that $g = 7$. By Claim~\ref{c:3}, $|S_X^*| \ge 5\delta - 11$. Thus, by Claim~\ref{c:1} and Claim~\ref{c:2}(a), $5\delta - 9 \ge |S| \ge |S_X| + |S_X^*| \ge 5\delta - 9$. Hence, we must have equality throughout this inequality chain, implying that $|S| = 5\delta - 9$, and so $|S_X| = 2$ and $|S_X^*| = 5\delta - 11$. Further, $S_X = \{x_1,x_3\}$. Since $x_4 \notin S_X$, the vertex $x_4$ is by Claim~\ref{c:4}(b) adjacent to $x_1$, $x_2$ or $x_3$. The girth condition implies in this case that $x_4$ has at most one neighbor in $X_{\le 3}$. This in turn implies that $m(D) \le 5$. Hence, by Claim~\ref{c:2}(c), $|S_X^*| \ge 6\delta - 10$, a contradiction. Hence, the girth $g = 8$.

The girth eight requirement implies that the vertex $x_i$ has at most one neighbor in $X_{\le 3}$ for every $i \in \{4,5,6\}$. This in turn implies that  $m(D) \le 7$. Hence, by Claim~\ref{c:2}(c), $|S_X^*| \ge 6\delta - 13$. Thus, by Claim~\ref{c:1} and Claim~\ref{c:2}(a), $5\delta - 11 \ge |S| \ge |S_X| + |S_X^*| \ge 6\delta - 11$. Hence, we must have equality throughout this inequality chain, implying that $|S| = 6\delta - 11$, and so $|S_X| = 2$ and $|S_X^*| = 6\delta - 13$. Further, $S_X = \{x_1,x_3\}$. Since $x_4 \notin S_X$, the vertex $x_4$ is adjacent to $x_1$, $x_2$ or $x_3$. The girth condition implies that $x_4$ is adjacent to exactly one vertex in $X_{\le 3}$. Since neither $x_5$ nor $x_6$ belongs to $S_X$, Claim~\ref{c:4}(b) implies that $x_i$ is adjacent to some vertex $x_j$ where $j \in [i-1]$ for $i \in \{5,6\}$. This, together with the girth eight requirement, implies that the vertex $x_5$ has exactly one neighbor in $X_{\le 4}$, and the vertex $x_6$ has at most two neighbors in $X_{\le 5}$. Therefore, $m(D) \le 6$. Hence, by Claim~\ref{c:2}(c), $|S_X^*| \ge 6\delta - 12$, a contradiction.~\smallqed

\medskip
By Claim~\ref{c:6.2}, the vertex $x_3$ has no neighbor in $X_{\le 2}$. If $x_4$ has at most one neighbor in $X_{\le 3}$, then $m(D) \le 5$. If $x_4$ has two neighbors in $X_{\le 3}$, then $x_4$ has one neighbor in $X_1 \cup X_2$ and a different neighbor in $X_3$. In this case, $x_5$ has at most one neighbor in each of $X_1 \cup X_2$ and $X_3 \cup X_4$, implying once again that $m(D) \le 5$. Thus if $g = 7$, then by Claim~\ref{c:2}(c), $|S_X^*| \ge 5\delta - 10$, and so $|S| \ge |S_X| + |S_X^*| \ge 5\delta -8$, contradicting Claim~\ref{c:1}(a). Therefore, $g = 8$.

\begin{subclaim}\label{c:6.3}
The vertex $x_4$ is adjacent to no vertex $x_j$ where $j \in [3]$.
\end{subclaim}
\proof
Suppose, to the contrary, that $x_4$ is adjacent to some vertex $x_j$, where $j \in [3]$. We show first that $x_4$ is adjacent to $x_3$.

\begin{unnumbered}{Claim~\ref{c:6.3}.1}
$x_4$ is adjacent to neither $x_1$ nor $x_2$.
\end{unnumbered}
\proof
Suppose, to the contrary, that $x_4$ is adjacent to $x_1$ or $x_2$. In this case, $x_i$ has at most one neighbor in $X_1 \cup X_2 \cup X_4$ for $i \in \{5,6\}$. Recall that $G_D$ is triangle-free. If $x_5$ has a neighbor in $X_3$, then $x_6$ has at most one neighbor in $X_3 \cup X_5$, by the girth condition, implying that $m(D) \le 7$. If $x_5$ has no neighbor in $X_3$, then $x_6$ has at most two neighbors in $X_3 \cup X_5$, implying once again that $m(D) \le 7$. Hence, by Claim~\ref{c:2}(c), $|S_X^*| \ge 6\delta - 13$. Thus, by Claim~\ref{c:1} and Claim~\ref{c:2}(a), $5\delta - 11 \ge |S| \ge |S_X| + |S_X^*| \ge 6\delta - 11$. Hence, we must have equality throughout this inequality chain, implying that $|S| = 6\delta - 11$, and so $|S_X| = 2$, $S_X = \{x_1,x_3\}$, and $|S_X^*| = 6\delta - 13$. Further, $m(D) = 7$, implying that $x_4$ has a neighbor in $X_{\le 3}$, and each of $x_5$ and $x_6$ has a neighbor in $X_1 \cup X_2 \cup X_4$. Further, at least one of $x_5$ and $x_6$ has a neighbor in $X_3$.

Since $x_5 \notin S_X$, the vertex $x_5$ is adjacent to some vertex $x_j$, where $j \in [4]$. If $x_5$ is adjacent to $x_3$, then, noting that $x_4$ has a neighbor in $X_{\le 3}$, the vertex $x_5$ and its two neighbors in $X_{\le 4}$ belong to a common cycle of length at most~$7$, a contradiction. Hence, $x_5$ is not adjacent to $x_3$ and is therefore adjacent to $x_1$, $x_2$ or $x_4$. The girth condition implies now that $x_5$ has no neighbor in $X_3$. Thus, $x_6$ has three neighbors in $X_{\le 3}$, one in each of $X_1 \cup X_2 \cup X_4$, $X_3$ and $X_5$. However, $x_6$ and its two neighbors in $X_{\le 5} \setminus X_3$ belong to a common cycle of length at most~$7$, a contradiction.~\smallqed

\medskip
By Claim~\ref{c:6.3}.1, $x_4$ is adjacent to $x_3$. By the girth condition, $x_i$ has at most one neighbor in each of $X_1 \cup X_2$ and $X_3 \cup X_4$ for $i \in \{5,6\}$, and so $x_i$ has at most two neighbors in $X_{\le 4}$. If $x_5$ has a neighbor in $X_1 \cup X_2$, then $x_6$ has at most one neighbor in $X_1 \cup X_2 \cup X_5$. If $x_5$ has a neighbor in $X_3 \cup X_4$, then $x_6$ has at most one neighbor in $X_3 \cup X_4 \cup X_5$. Hence, if $x_5$ has at least one neighbor in $X_{\le 4}$, then $x_6$ has at most two neighbors in $X_{\le 5}$. This implies that $m(D) \le 7$. Analogously as in the proof of Claim~\ref{c:6.3}.1, we deduce that $|S| = 6\delta - 11$, $|S_X^*| = 6\delta - 13$ and $S_X = \{x_1,x_3\}$. Further, $m(D) = 7$, implying that $x_4$ has a common neighbor with either $x_1$ or $x_2$, and $x_5$ has two neighbors in $X_{\le 4}$. Since $x_5 \notin S_X$, the vertex $x_5$ is adjacent to some vertex $x_j$, where $j \in [4]$. The above properties of the graph $G$ imply that the vertex $x_5$ and its two neighbors in $X_{\le 4}$ belong to a common cycle of length at most~$7$, a contradiction.  This completes the proof of Claim~\ref{c:6.3}.~\smallqed

\medskip
By Claim~\ref{c:6.3}, the vertex $x_4$ is not adjacent to $x_1$, $x_2$ or $x_3$. Hence by Claim~\ref{c:4}(b), $x_4 \in S_X$, and so $|S_X| \ge 3$ and $\{x_1,x_3,x_4\} \subseteq S_X$.

\begin{subclaim}\label{c:6.4}
The vertex $x_4$ has no neighbor in $X_{\le 3}$.
\end{subclaim}
\proof
Suppose, to the contrary, that $x_4$ has a neighbor in $X_{\le 3}$. If $x_4$ has a neighbor in $X_1 \cup X_2$, then analogous arguments as in the proof of Claim~\ref{c:6.3}.1 show that $m(D) \le 7$. Hence, by Claim~\ref{c:2}(c), $|S_X^*| \ge 6\delta - 13$. Thus, by Claim~\ref{c:2}(a), $|S| \ge |S_X| + |S_X^*| \ge 6\delta - 10$, a contradiction. Hence, $x_4$ has no neighbor in $X_1 \cup X_2$, and therefore $x_4$ has a neighbor in $X_3$. Thus, $x_5$ and $x_6$ have at most one neighbor in each of $X_1 \cup X_2$ and $X_3 \cup X_4$, implying that $m(D) \le 7$ and therefore, as before, that $|S_X^*| \ge 6\delta - 13$ and $|S| \ge 6\delta - 10$, a contradiction.~\smallqed

\medskip
We now return to the proof of Claim~\ref{c:6} one final time. By Claim~\ref{c:6.4}, the vertex $x_4$ has no neighbor in $X_{\le 3}$. Thus, by our earlier observation, the vertex $x_4$ is not adjacent to $x_1$, $x_2$ or $x_3$ in the graph $G_D$. As observed earlier, $x_1$ and $x_2$ are adjacent in $G_D$, and $x_3$ is not adjacent to $x_1$ or $x_2$ in $G_D$. Since $G_D$ is triangle-free, this implies that $m(D) \le 7$. Hence, by Claim~\ref{c:2}(c), $|S_X^*| \ge 6\delta - 13$. Thus, by Claim~\ref{c:2}(a), $|S| \ge |S_X| + |S_X^*| \ge 6\delta - 10$, a contradiction. This completes the proof of Claim~\ref{c:6}.~\smallqed

\medskip
By Claim~\ref{c:6}, the vertices $x_1$ and $x_2$ are not adjacent. Hence by Claim~\ref{c:4}(b), $x_2 \in S_X$, implying that $\{x_1,x_2\} \subseteq S_X$ and $|S_X| \ge 2$.

\begin{claim}\label{c:7}
The vertex $x_3$ is not adjacent to $x_1$ or $x_2$.
\end{claim}
\proof  Suppose, to the contrary, that $x_3$ is adjacent to $x_1$ or $x_2$.

\begin{subclaim}\label{c:7.1}
The vertex $x_3$ is adjacent to exactly one of $x_1$ and $x_2$.
\end{subclaim}
\proof Suppose, to the contrary, that $x_3$ is adjacent to both $x_1$ and $x_2$. In this case, $x_2$ has no neighbor in $X_1$. Further, $x_i$ has at most one neighbor in $X_{\le 3}$ for $i \in \{4,\ldots,g-2\}$. If $g = 7$, this implies that $m(D) \le 5$, and so, by Claim~\ref{c:2}(c), $|S_X^*| \ge 6\delta - 10$. Thus, by Claim~\ref{c:2}(a), $|S| \ge |S_X| + |S_X^*| \ge 6\delta - 8$, a contradiction. Hence, $g = 8$, implying by the triangle-freeness of $G_D$ and by our earlier observations that $m(D) \le 7$. Thus, by Claim~\ref{c:2}(c), $|S_X^*| \ge 6\delta - 13$. If $x_4$ is adjacent to some vertex $x_j$, where $j \in [3]$, then $x_i$ has at most one neighbor in $X_{\le 4}$ for $i \in \{5,6\}$, implying that $m(D) \le 6$, and so, by Claim~\ref{c:2}(c), $|S_X^*| \ge 6\delta - 12$ and therefore, by Claim~\ref{c:2}(a), $|S| \ge |S_X| + |S_X^*| \ge 6\delta - 10$, a contradiction. Hence, $x_4$ is not adjacent to some vertex $x_j$, where $j \in [3]$. Thus, by Claim~\ref{c:4}(b), $x_4 \in S_X$, and so $|S_X| \ge 3$. Therefore, by Claim~\ref{c:2}(a), $|S| \ge |S_X| + |S_X^*| \ge 3 + (6\delta - 13) = 6\delta - 10$, a contradiction.~\smallqed

\medskip
By Claim~\ref{c:7.1}, the vertex $x_3$ is adjacent to exactly one of $x_1$ and $x_2$. Let $\{x_1,x_2\} = \{x_{i_1},x_{i_2}\}$, where $x_3$ is adjacent to $x_{i_1}$. The girth requirement implies that each of $x_{i_2}$, $x_4$, $x_5$ and, if $g = 8$, $x_6$ has at most one neighbor in $X_{i_1} \cup X_3$.

\begin{subclaim}\label{c:7.2}
$g = 8$.
\end{subclaim}
\proof
Suppose, to the contrary, that $g = 7$.  As observed earlier, each of $x_{i_2}$, $x_4$ and $x_5$ has at most one neighbor in $X_{i_1} \cup X_3$. Therefore since $G_D$ is triangle-free, $m(D) \le 6$. Hence, by Claim~\ref{c:2}(c), $|S_X^*| \ge 5\delta - 11$. Thus, by Claim~\ref{c:1} and Claim~\ref{c:2}(a), $5\delta - 9 \ge |S| \ge |S_X| + |S_X^*| \ge 5\delta - 9$. Hence, we must have equality throughout this inequality chain, implying that $|S| = 5\delta - 9$, and so $|S_X| = 2$, $S_X = \{x_1,x_2\}$, and $|S_X^*| = 5\delta - 11$. Further, $m(D) = 6$, implying that each of $x_{i_2}$, $x_4$ and $x_5$ has a neighbor in $X_{i_1} \cup X_3$. Since $x_4 \notin S_X$, the vertex $x_4$ is adjacent to some vertex $x_j$, where $j \in [3]$. If $x_4$ is adjacent to $x_{i_1}$ or $x_3$, then $x_4$ has no neighbor in $X_{i_2}$, and $x_5$ has at most one neighbor in $X_{i_1} \cup X_3 \cup X_4$, implying that $m(D) \le 5$, a contradiction. Hence, $x_4$ is adjacent to $x_{i_2}$. However, since both $x_{i_2}$ and $x_4$ have a neighbor in $X_{i_1} \cup X_3$, we contradict the girth condition.~\smallqed

\medskip
By Claim~\ref{c:7.2}, $g = 8$.

\begin{subclaim}\label{c:7.3}
There is no edge in $G_D$ joining $x_{i_2}$ to $x_{i_1}$ or $x_3$.
\end{subclaim}
\proof
Suppose, to the contrary, that $x_{i_2}$ is adjacent to $x_{i_1}$ or $x_3$ in $G_D$. This implies that
$x_1$ and $x_2$ have a common neighbor or $x_3$ has a common neighbor with $x_{i_2}$. Thus, $x_i$ has at most one neighbor in $X_{\le 3}$ for $i \in \{4,5,6\}$. Therefore since $G_D$ is triangle-free, $m(D) \le 7$. Hence, by Claim~\ref{c:2}(c), $|S_X^*| \ge 6\delta - 13$. Thus, by Claim~\ref{c:1} and Claim~\ref{c:2}(a), $6\delta - 11 \ge |S| \ge |S_X| + |S_X^*| \ge 6\delta - 11$. Hence, we must have equality throughout this inequality chain, implying that $|S| = 6\delta - 11$, and so $|S_X| = 2$, $S_X = \{x_1,x_2\}$, and $|S_X^*| = 5\delta - 13$. Further, $m(D) = 7$, implying that each of $x_4$, $x_5$ and $x_6$ has at most one neighbor in $X_{\le 3}$. Since $x_4 \notin S_X$, the vertex $x_4$ is adjacent to some vertex $x_j$, where $j \in [3]$. If $x_4$ is adjacent to $x_{i_1}$ or $x_3$, then $x_4$ has no neighbor in $X_{i_2}$, and both $x_5$ and $x_6$ have at most one neighbor in $X_{\le 4}$, implying that $m(D) \le 6$, a contradiction. Hence, $x_4$ is adjacent to $x_{i_2}$, and, by the girth condition, has no neighbor in $X_{i_1} \cup X_3$. Since $x_5 \notin S_X$, the vertex $x_5$ is adjacent to some vertex $x_j$, where $j \in [4]$, implying that $x_5$ has at exactly one neighbor in $X_{\le 4}$ and $x_6$ has at most two neighbors in $X_{\le 5}$. Once again, this implies that $m(D) \le 6$, a contradiction.~\smallqed

\medskip
By Claim~\ref{c:7.3}, there is no edge in $G_D$ joining $x_{i_2}$ to $x_{i_1}$ or $x_3$. Thus, the vertices $x_2$ has no neighbor in $X_1$, and the vertex $x_3$ has exactly one neighbor in $X_{\le 2}$, namely the vertex $x_{i_1}$.

\begin{subclaim}\label{c:7.4}
The vertex $x_4$ is adjacent to no vertex $x_j$ where $j \in [3]$.
\end{subclaim}
\proof
Suppose, to the contrary, that $x_4$ is adjacent to some vertex $x_j$, where $j \in [3]$. In this case, $x_i$ has at most two neighbors in $X_{\le 4}$ for $i \in \{5,6\}$. Further, if $x_5$ has two neighbors in $X_{\le 4}$, then $x_6$ has at most two neighbors in $X_{\le 5}$, implying that $m(D) \le 7$. If $x_5$ has at most one neighbor in $X_{\le 4}$, then $x_6$ has at most three neighbors in $X_{\le 5}$, implying once again $m(D) \le 7$. Hence, $m(D) \le 7$. Thus, by Claim~\ref{c:2}(c), $|S_X^*| \ge 5\delta - 13$. Analogous arguments as before (see, for example, the proof of Claim~\ref{c:7.3}) imply that $|S| = 6\delta - 11$, $S_X = \{x_1,x_2\}$, and $|S_X^*| = 5\delta - 13$. Further, $m(D) = 7$, implying in particular that $x_4$ has two neighbors in $X_{\le 3}$. Since $x_5 \notin S_X$, the vertex $x_5$ is adjacent to some vertex $x_j$, where $j \in [4]$, implying that $x_5$ has at exactly one neighbor in $X_{\le 4}$ and $x_6$ has at most two neighbors in $X_{\le 5}$. This implies that $m(D) \le 6$, a contradiction.~\smallqed

\medskip
By Claim~\ref{c:7.4}, the vertex $x_4$ is adjacent to no vertex $x_j$ where $j \in [3]$. Hence by Claim~\ref{c:4}(b), $x_4 \in S_X$, and so $\{x_1,x_2,x_4\} \subseteq S_X$. Thus, $|S_X| \ge 3$.
Suppose that $x_4$ has at least one neighbor in $X_{i_1} \cup X_3$. The girth condition now implies that each of $x_{i_2}$, $x_5$ and $x_6$ has at most one neighbor in $X_{i_1} \cup X_3 \cup X_4$. Further, since $G_D$ is triangle-free, at most two of the edges $x_{i_2}x_5$, $x_{i_2}x_6$, and $x_5x_6$ are present in $G_D$. Thus, $m(D) \le 7$. Hence, by Claim~\ref{c:2}(c), $|S_X^*| \ge 5\delta - 13$. Thus, by Claim~\ref{c:2}(a), $|S| \ge |S_X| + |S_X^*| \ge 3 + (6\delta - 13) \ge 6\delta - 10$, a contradiction. Therefore, $x_4$ has no neighbor in $X_{i_1} \cup X_3$. Hence, in $G_D$ there is no edge joining a vertex in $\{x_{i_2},x_4\}$ and a vertex in $\{x_{i_1},x_3\}$. As observed earlier, each of $x_5$ and $x_6$ has at most one neighbor in $X_{i_1} \cup X_3$, and therefore there are at most two edges in $G_D$ joining vertices in $\{x_5,x_6\}$ and vertices in $\{x_{i_1},x_3\}$. Further since $G_D$ is triangle-free, there are at most four edges in the subgraph of $G_D$ induced by $\{x_{i_2},x_4,x_5,x_6\}$. Therefore, $m(D) \le 7$, once again producing a contradiction. This completes the proof of Claim~\ref{c:7}.~\smallqed

\medskip
By Claim~\ref{c:7}, the vertex $x_3$ is not adjacent to $x_1$ or $x_2$. By Claim~\ref{c:4}(b), $x_3 \in S_X$, implying that $\{x_1,x_2,x_3\} \subseteq S_X$ and $|S_X| \ge 3$.

\begin{claim}\label{c:8}
The following holds. \\[-28pt]
\begin{enumerate}
\item $g = 8$.
\item $S_X = \{x_1,x_2,x_3\}$.
\item $m(D) = 8$.
\end{enumerate}
\end{claim}
\proof If $g = 7$, then by Claim~\ref{c:3}, $|S_X^*| \ge 5\delta - 11$, and so by Claim~\ref{c:2}(a), $|S| \ge |S_X^*| + |S_X| \ge 5\delta - 8$, a contradiction. Therefore, $g = 8$. This proves Part~(a). We prove next that $|S_X^*| \ge 6\delta - 14$. Suppose, to the contrary, that $|S_X^*| < 6\delta - 14$. By Claim~\ref{c:3}(b), this implies that $|S_X^*| = 6\delta - 15$. This in turn implies that the triangle-free graph $G_D$ of order~$6$ has maximum possible size, namely $9$, and is therefore the graph $K_{3,3}$. By Claim~\ref{c:1} and Claim~\ref{c:2}(a), $6\delta - 11 \ge |S| \ge |S_X| + |S_X^*| = |S_X| + 6\delta - 15$, implying that $|S_X| \le 4$. Since $\{x_1,x_2,x_3\} \subseteq S_X$, at least two vertices in $\{x_4,x_5,x_6\}$ do not belong to the set $S_X$. Let $x_{i_0}$ be such a vertex in $\{x_4,x_5,x_6\}$. Thus, $x_{i_0}$ is adjacent in $G$ to some vertex $x_{i_1}$, where $i_1 < i_0$. Let $x_{i_0}x_{i_1}x_{i_2}x_{i_3}x_{i_0}$ be a $4$-cycle in $G_D \cong K_{3,3}$ that contains the edge $x_{i_0}x_{i_1}$. We note that $x_{i_j}$ and $x_{i_{j+1}}$ are either adjacent in $G$ or have a common neighbor in $G$ for $j \in \{1,2,3\}$, where addition is taken modulo~$4$. Thus, since $x_{i_0}x_{i_1}$ is an edge in $G$, this implies that $x_{i_0}$ and $x_{i_1}$ belong to a cycle of length at most~$7$ in $G$, a contradiction. Hence, $|S_X^*| \ge 6\delta - 14$. Thus, by Claim~\ref{c:1} and Claim~\ref{c:2}(a), $6\delta - 11 \ge |S| \ge |S_X| + |S_X^*| \ge 3 + (6\delta - 14) = 6\delta - 11$. Hence, we must have equality throughout this inequality chain, implying that $|S| = 6\delta - 11$, and so $|S_X| = 3$, $S_X = \{x_1,x_2,x_3\}$. Further, $|S_X^*| = 5\delta - 14$, implying that $m(D) = 8$. This proves Parts~(b) and~(c).~\smallqed

\begin{claim}\label{c:9}
The vertex $x_4$ is adjacent to exactly one of $x_1$, $x_2$ and $x_3$.
\end{claim}
\proof By Claim~\ref{c:8}(b), the vertex $x_4 \notin S_X$, and so $x_4$ is adjacent to $x_1$, $x_2$ or $x_3$. Suppose, to the contrary, that $x_4$ is adjacent to at least two of $x_1$, $x_2$ and $x_3$. Let $\{x_1,x_2,x_3\} = \{x_{i_1},x_{i_2},x_{i_3}\}$, where $x_4$ is adjacent to $x_{i_1}$ and $x_{i_2}$. The girth requirement implies that each of $x_{i_3}$, $x_5$ and $x_6$ has at most one neighbor in $X_{i_1} \cup X_{i_2} \cup X_4$. Further, since $G_D$ is triangle-free, at most two of the edges $x_{i_3}x_5$, $x_{i_3}x_6$, and $x_5x_6$ are present in $G_D$. Thus, $m(D) \le 7$, contradicting Claim~\ref{c:8}(c).~\smallqed

\medskip
By Claim~\ref{c:9}, the vertex $x_4$ is adjacent to exactly one of $x_1$, $x_2$ and $x_3$. Let $\{x_1,x_2,x_3\} = \{x_{i_1},x_{i_2},x_{i_3}\}$, where $x_4$ is adjacent to $x_{i_1}$.
Suppose that $x_{i_2}$ has a common neighbor with $x_{i_1}$ or $x_4$ in $G$. The girth requirement implies that each of $x_{i_3}$, $x_5$ and $x_6$ has at most one neighbor in $X_{i_1} \cup X_{i_2} \cup X_4$. Further, since $G_D$ is triangle-free, at most two of the edges $x_{i_3}x_5$, $x_{i_3}x_6$, and $x_5x_6$ are present in $G_D$. Thus, $m(D) \le 7$, contradicting Claim~\ref{c:8}(c). Hence, $x_{i_2}$ has no common neighbor with $x_{i_1}$ or $x_4$ in $G$. Analogously, $x_{i_3}$ has no common neighbor with $x_{i_1}$ or $x_4$ in $G$. Hence in $G_D$ there is no edge joining a vertex in $\{x_{i_1},x_4\}$ and a vertex in $\{x_{i_2},x_3\}$. Each of $x_5$ and $x_6$ has at most one neighbor in $X_{i_1} \cup X_4$, and therefore there are at most two edges in $G_D$ joining vertices in $\{x_5,x_6\}$ and vertices in $\{x_{i_1},x_4\}$. Further since $G_D$ is triangle-free, there are at most four edges in the subgraph of $G_D$ induced by $\{x_{i_2},x_{i_3},x_5,x_6\}$. Therefore, $m(D) \le 7$, contradicting Claim~\ref{c:8}(c). This completes the proof of Theorem~\ref{thm:main1a}.~\qed


\medskip
\begin{thebibliography}{99}

\bibitem{AIM-Workshop}
AIM Special Work Group, Zero forcing sets and the minimum rank of graphs. \textit{Linear Algebra Appl.} \textbf{428}(7) (2008), 1628--1648.


\bibitem{k-Forcing}
D. Amos, Y. Caro, R. Davila, and R. Pepper, Upper bounds on the $k$-forcing number of a graph. \textit{Discrete Applied Math.} \textbf{181} (2015), 1--10.

\bibitem{Barioli13}
F. Barioli,  W. Barrett, S. M. Fallat, T. Hall, L. Hogben, B. Shader, P. van den Driessche, and H. van der Holst, Parameters related to tree-width, zero forcing, and maximum nullity of a graph. \textit{J. Graph Theory} \textbf{72}(2) (2013), 146--177.


\bibitem{quantum1}
D. Burgarth and V. Giovannetti, Full control by locally induced relaxation. \textit{Physical Review Letters} \textbf{99}(10) (2007), 100501.



\bibitem{logic1}
D. Burgarth, V. Giovannetti,  L. Hogben,  S. Severini,  and M. Young.
\newblock Logic circuits from zero forcing.
\newblock {arXiv preprint arXiv:1106.4403}, 2011.


\bibitem{Dynamic Forcing}
Y. Caro and R. Pepper.
\newblock Dynamic approach to k-forcing.
\newblock {\em Theory and Applications of Graphs}, Volume 2: Iss. 2, Article 2, 2015.

\bibitem{ChSu05} S. Chandran and C. Subramanian, Girth and treewidth. \textit{J. Combinatorial Theory B} \textbf{93} (2005), 23--32.




\bibitem{zf_np}
C. Chekuri and N. Korula.
\newblock A graph reduction step preserving element-connectivity and applications.
\newblock {\em Automata, Languages and Programming}, 254--265. Springer 2009.


\bibitem{Davila Thesis}
R. Davila.
\newblock Bounding the forcing number of a graph.
\newblock{\em Rice University Masters Thesis}, 2015.


\bibitem{DaHeMaPe16}
R. Davila, M. A. Henning, C. Magnant, and R. Pepper, Bounds on the connected forcing number of a graph, manuscript.



\bibitem{Davila Kenter}
R. Davila, and F. Kenter.
\newblock Bounds for the zero forcing number of a graph with large girth.
\newblock {\em Theory and Applications of Graphs}, Volume 2, Issue 2, Article 1, 2015.

\bibitem{GaKwLa93} D. K. Garnick, Y. H. Harris Kwong, and F. Lazebnik, Extremal graphs without three-cycles or four-cycles. \textit{J. Graph Theory} \textbf{17}(5) (1993), 633--645.


\bibitem{Genter1} M. Gentner, L. D. Penso, D. Rautenbach, and U. S. Souzab, Extremal values and bounds for the zero forcing number. \textit{Discrete Applies Math.} \textbf{214} (2016), 196--200.


\bibitem{Genter2} M. Gentner and D. Rautenbach, Some bounds on the zero forcing number of a graph, manuscript. https://arxiv.org/pdf/1608.00747.pdf



\bibitem{Edholm}
C. Edholm, L. Hogben, J. LaGrange, and D. Row, Vertex and edge spread of zero forcing number, maximum nullity, and minimum rank of a graph. \textit{Linear Algebra Appl.} \textbf{436}(12) (2012), 4352--4372.



\bibitem{powerdom3}
T. W. Haynes,  S. T. Hedetniemi, S. T. Hedetniemi, and M. A. Henning, Domination in graphs applied to electric power networks. \textit{SIAM J. Discrete Math.} \textbf{15}(4) (2002), 519--529.



\bibitem{MHAYbookTD} M. A. Henning and A. Yeo, \emph{Total domination in graphs (Springer Monographs in Mathematics)}.  ISBN-13: 978-1461465249, 2013.



\bibitem{zf_np2}
M. Trefois and J.C. Delvenne.
\newblock Zero forcing sets, constrained matchings and minimum rank.
\newblock {\em Linear and Multilinear Algebra}, 2013.


\bibitem{powerdom2}
M. Zhao,  L. Kang, and G. Chang, Power domination in graphs. \textit{Discrete Math.} \textbf{306} (2006), 1812--1816.




\end{thebibliography}
\end{document}